\documentclass[12pt, a4paper]{article}
\usepackage{amsmath,amssymb,amsbsy,amsfonts,amsthm,latexsym,
                          amsopn,amstext,amsxtra,euscript,amscd}

\begin{document}

\newtheorem{lem}{Lemma}
\newtheorem{lemma}[lem]{Lemma}
\newtheorem{prop}{Proposition}
\newtheorem{thm}{Theorem}
\newtheorem{theorem}[thm]{Theorem}

\title{On multiplicative congruences}

\author{{M. Z. Garaev}
\\
\normalsize{Instituto de Matem{\'a}ticas}
\\
\normalsize{Universidad Nacional Aut\'onoma de M\'exico}
\\
\normalsize{Campus Morelia, Apartado Postal 61-3 (Xangari)}
\\
\normalsize{C.P. 58089, Morelia, Michoac{\'a}n, M{\'e}xico} \\
\normalsize{\tt garaev@matmor.unam.mx}\\
}

\date{}

\pagenumbering{arabic}

\maketitle

\begin{abstract}
Let $\varepsilon$ be a fixed positive quantity, $m$ be a large
integer, $x_j$ denote integer variables. We prove that for any
positive integers $N_1,N_2,N_3$ with $N_1N_2N_3>m^{1+\varepsilon}, $
the set
$$
\{\, x_1x_2x_3 \pmod m: \quad x_j\in [1,N_j]\,\}
$$
contains almost all the residue classes modulo $m$ (i.e., its
cardinality is equal to $m+o(m)$). We further show that if $m$ is
cubefree, then for any positive integers $N_1,N_2,N_3,N_4$ with $
N_1N_2N_3N_4>m^{1+\varepsilon}, $ the set
$$
\{\, x_1x_2x_3x_4 \pmod m: \quad x_j\in [1,N_j]\,\}
$$
also contains almost all the residue classes modulo $m.$

Let $p$ be a large prime parameter and let
$p>N>p^{63/76+\varepsilon}.$ We prove that for any nonzero integer
constant $k$ and any integer $\lambda\not\equiv 0\pmod p$ the
congruence
$$
p_1p_2(p_3+k)\equiv \lambda\pmod p
$$
admits $(1+o(1))\pi(N)^3/p$ solutions in prime numbers $p_1, p_2,
p_3\le N.$
\end{abstract}

\paragraph*{2000 Mathematics Subject Classification:} 11L40

\newpage

\section {Introduction}

In our works~\cite{Gar1,Gar2} we applied large value results of
character sums to a concrete multiplicative ternary congruence and
by this mean improved one of the results of Friedlander and
Shparlinski~\cite{FrSh}. In the present paper we examine those
arguments in application to some other multiplicative congruences.

Everywhere below $\varepsilon$ denotes a small fixed positive
quantity, $m$ is a large integer parameter.

\begin{theorem}
\label{thm:main1} Let $N_1,N_2, N_3$ be positive integers with
$$
N_1N_2N_3>m^{1+\varepsilon}.
$$
Then for some $\delta=\delta(\varepsilon)>0$ we have
$$
\#\{\,x_1x_2x_3 \pmod m: \quad x_j\in [1, N_j]\,\} \, =\,
m+O(m^{1-\delta}).
$$
\end{theorem}

In the statement of Theorem~\ref{thm:main1} the condition $N_1N_2
N_3>m^{1+\varepsilon}$ can not be relaxed to $N_1N_2N_3>Cm,$ no
matter how large the constant $C$ is. We also note that if $m=qn,$
where $q$ is a prime number approximately several times bigger than
$n^{(1+\varepsilon)/(2-\varepsilon)},$ then
$q>m^{(1+\varepsilon)/3}$ and hence non of the $n$ numbers
$q,2q,\ldots, nq$ can be represented in the form $x_1x_2x_3\pmod m$
with $x_j\le m^{(1+\varepsilon)/3}.$ In particular the exponent of
$m$ inside of the $O$-symbol can not be replaced by a constant
smaller than $(2-\varepsilon)/3.$

It is known~\cite{GK} that the set
$$
{\cal A}_2:=\{\,x_1x_2 \pmod m: \quad x_1,x_2\in [1,
m^{1/2+\varepsilon}]\,\}
$$
contains almost all the elements of the residue ring $\mathbb{Z}_m.$
Theorem~\ref{thm:main1} implies that the set
$$
{\cal A}_3:=\{\,x_1x_2x_3 \pmod m: \quad x_1,x_2,x_3\in [1,
m^{1/3+\varepsilon}]\,\}
$$
also contains almost all the elements of $\mathbb{Z}_m.$ Another
consequence of Theorem~\ref{thm:main1} is that for any sufficiently
large integer $m,$ any invertible element $\lambda\in
\mathbb{Z}^*_m$ can be represented in the form
$$
\lambda\equiv x_1x_2x_3x_4x_5x_6\pmod m
$$
for some positive integers $x_1,x_2,\ldots,x_6$ with
$\max\limits_{1\le i\le 6}x_i\le m^{1/3+\varepsilon}.$

\begin{theorem}
\label{thm:main11} Let $m$ be cubefree, $N_1,N_2,N_3, N_4$ be
positive integers with
$$
N_1N_2N_3N_4>m^{1+\varepsilon}.
$$
Then for some $\delta=\delta(\varepsilon)>0$ we have
$$
\#\{\,x_1x_2x_3x_4 \pmod m: \quad x_j\in [1, N_j]\,\} \, = \,
m+O(m^{1-\delta}).
$$
\end{theorem}
In particular, for cubefree $m$ the set
$$
{\cal A}_4:=\{\,x_1x_2x_3x_4 \pmod m: \quad  x_1,x_2,x_3,x_4\in [1,
m^{1/4+\varepsilon}]\,\}
$$
contains almost all the  elements of $\mathbb{Z}_m$ and also
contains almost all the elements of $\mathbb{Z}^*_m.$ This implies
that, for any sufficiently large integer $m,$ any element
$\lambda\in \mathbb{Z}_m^*$ is representable in the form
$$
\lambda\equiv x_1x_2x_3x_4x_5x_6x_7x_8\pmod m
$$
for some positive integers $x_1,x_2,\ldots,x_8$ with
$\max\limits_{1\le i\le 8}x_i\le m^{1/4+\varepsilon}.$ It would be
interesting to reduce the number of variables in the latter
statement to $7.$

\begin{theorem}
\label{thm:main2} Let $p$ be a large prime parameter, $k$ be a
nonzero integer constant and $\lambda$ be an integer coprime to $p.$
If $p>N>p^{63/76+\varepsilon}$ then the congruence
$$
p_1p_2(p_3+k)\equiv \lambda\pmod p
$$
has $(1+o(1))\pi(N)^3/p$ solutions in primes $p_1, p_2, p_3\le N.$
\end{theorem}

Theorem~\ref{thm:main2} quantitatively complements Theorem 6 from
the work of Friedlander, Kurlberg and Shparlinski~\cite{FKSh}.

In what follows, the letters
$\varepsilon',\varepsilon'',\varepsilon''',\varepsilon_1$ are used
to denote some positive fixed quantities chosen in obvious ways. The
letters $x_j, y_j, u_j, t$ denote integer numbers.

\section{Character sum estimates}

In the proofs of Theorems~\ref{thm:main1},~\ref{thm:main11} we will
use well-known character sum estimates of Burgess~\cite{Bur1, Bur2}:
if $N>m^{1/3+\varepsilon}$  then there exists
$\delta=\delta(\varepsilon)>0$ such that for any nonprincipal
character $\chi\pmod m$ we have
$$
\left|\sum_{x\le N}\chi(x)\right|\le Nm^{-\delta}.
$$
In the case when $m$ is cubefree, the condition
$N>m^{1/3+\varepsilon}$ can be relaxed to $N>m^{1/4+\varepsilon}.$

To prove Theorem~\ref{thm:main2} we shall use Vinogradov's bound on
character sums over shifted primes. Let $k$ be a fixed nonzero
integer constant, $\chi$ be a nonprincipal character modulo $p.$
Then Vinogradov's work~\cite{Vin} implies that in the range $1\le
N<p$ one has
\begin{equation}
\label{eqn:Vin} \left|\sum_{p'\le N}\chi(p'+k)\right|\lessapprox
p^{1/4}N^{2/3},
\end{equation}
where $p'$ denotes prime numbers. Here and below, we use the
notation $L\lessapprox M$ to indicate that for any fixed
$\varepsilon>0$ there exists a constant $c=c(\varepsilon)$ such that
$L\le cMp^{\varepsilon}$ (or $L\le cMm^{\varepsilon}$ in the proofs
of Theorems~\ref{thm:main1},~\ref{thm:main11}).

The bound~\eqref{eqn:Vin} is nontrivial when
$N>p^{3/4+\varepsilon}.$ We mention that Karatsuba's work~\cite{Ka}
implies a nontrivial bound in the wider range
$N>p^{1/2+\varepsilon}.$ Since we deal with larger values of $N,$
the estimate~\eqref{eqn:Vin} will be more profitable.

\section{Large values of character sums}

Having character sum estimates under hands, one can apply
Karatsuba's  method from~\cite{Ka2} to  derive a variety of results
on solvability of multiplicative ternary congruences and find
asymptotic formulas for the number of their solutions. Our theorems,
however, can not be obtained from the direct application of
Karatsuba's method combined with Burgess' and Vinogradov's character
sum estimates. One main ingredient in our proofs is Huxley's
refinement of the Hal\'asz-Montgomery method for large value results
of Dirichlet polynomials. Our present application of this theory can
be compared with Lemma~4 of Friedlander and Iwaniec~\cite{FrIw}. For
our purposes it suffices the following simplest form of it. Let
$a_n$ be numbers with $|a_n|\lessapprox 1,$ let $0< V\le N$ and let
$R$ be the number of characters $\chi\pmod p$ for which
$$
\Bigl|\sum_{n=N+1}^{2N}a_n\chi(n)\Bigr|\ge V.
$$
Then Huxley's refinement implies that
\begin{equation}
\label{eqn:Hux1} R\lessapprox \frac{N^2}{V^2}+\frac{pN^4}{V^6},
\end{equation}
see Mongomery~\cite{Mont}, Huxley~\cite{Hux}, Huxley and
Jutila~\cite{HJ}, Jutila~\cite{J}.

The estimate~\eqref{eqn:Hux1} will be used in the proof of
Theorem~\ref{thm:main2}. A suitable version of it can also be used
to prove Theorems~\ref{thm:main1},~\ref{thm:main11}, but in this
case we can present the proof in a relatively more elementary
language, so that it will be more self-contained.

\section{Proof of Theorem~\ref{thm:main1}}

It suffices to prove the following lemma:

\begin{lemma}
\label{lem:main1} Let $N_1N_2N_3>m^{1+\varepsilon}.$ Then there are
only $O(m^{1-\varepsilon_1})$ elements  $\lambda\in \mathbb{Z}_m^*$
such that
$$
\lambda\not \in\{x_1x_2x_3\pmod m: \quad x_j\in [1,N_j]\}.
$$
\end{lemma}

Indeed, assume that Lemma~\ref{lem:main1} is proved and we show how
to derive Theorem~\ref{thm:main1} from this lemma.

In the condition of Theorem~\ref{lem:main1} we can assume that
$N_1>m^{(1+\varepsilon)/3}.$ Denote by ${\cal H}$ the set of all
elements $\lambda\in \mathbb{Z}_m$ such that
$$
\lambda\not \in\{x_1x_2x_3\pmod m: \quad x_j\in [1,N_j]\}.
$$
For a given divisor $d|m,$ let ${\cal H}_d$ be the set of all
elements of ${\cal H}$ such that $(h,m)=d$ for any $h\in {\cal
H}_d.$ Since $|{\cal H}_d| \le m/d,$ we have
\begin{equation}
\label{eqn: H}
|{\cal H}|=\sum_{d|m}|{\cal H}_d|=\sum_{\substack{d|m\\
d<m^{\varepsilon'}}}|{\cal H}_d|+\sum_{\substack{d|m\\
d\ge m^{\varepsilon'}}}m/d=\sum_{\substack{d|m\\
d<m^{\varepsilon'}}}|{\cal H}_d|+O(m^{1-0.5\varepsilon'}),
\end{equation}
where $\varepsilon'=0.1\varepsilon$ say. We estimate ${\cal H}_d$
for $d<m^{\varepsilon'}.$ By the definition,
$$
{\cal H}_d=\{\,dx\pmod m: \quad x\in {\cal B}_d\,\} \qquad {\rm for
\,\,\, some} \quad {\cal B}_d\subset Z_{m/d}^*.
$$
Since
$$
{\cal H}_d\cap \{x_1x_2x_3\pmod m: \quad x_j\in [1,N_j]\}=\emptyset,
$$
taking $x_1=dy_1,\, y_1\in [1, N_1/d]$ we get that
$$
{\cal  B}_d\cap \{x_1x_2x_3\pmod {(m/d)}: \, x_1\in [1,N_1/d],\,
x_2\in [1,N_2], \, x_3\in [1,N_3]\}=\emptyset.
$$
Since ${\cal  B}_d\subset Z_{m/d}^*$ and
$(N_1/d)N_2N_3>(m/d)^{1+0.5\varepsilon},$ we can apply
Lemma~\ref{lem:main1} with $m$ replaced by $m/d$ and $N_1$ replaced
by $[N_1/d]$ and deduce that
$$
|{\cal  H}_d|=|{\cal  B}_d|=O(m^{1-\varepsilon''}),\quad
\varepsilon''>0.
$$
Incorporating this into~\eqref{eqn: H}, we conclude that
$$
|{\cal  H}|=O(m^{1-\varepsilon'''}), \quad \varepsilon'''>0.
$$

Thus, it suffices to prove Lemma~\ref{lem:main1}. We can assume that
$m^{0.1\varepsilon}<N_j<m$ for all $j.$ Indeed, if say
$N_1<m^{0.1\varepsilon},$ then $N_2N_3>m^{1+0.9\varepsilon}$ and we
simply can take $x_1=1$ and look for $x_2=y_1y_2$ with $y_1,y_2\in
[1, N_2^{1/2}].$

Substituting $x_1\to ux_1$ and manipulating with $\varepsilon$ it
suffices to show that if $N_1N_2N_3>m^{1+\varepsilon}$ then for some
$\varepsilon_1>0,$
$$
\#\{ux_1x_2x_3\pmod m: \, u\in [1, U], \,x_j\in [1,
N_j]\}=m+O(m^{1-\epsilon_1}),
$$
where  $U=[m^{1/n}], \, n=[10/\varepsilon].$  We can assume that
$N_1>m^{(1+\varepsilon)/3}.$ From the Burgess character sum
estimate, there exists a positive quantity
$\delta=\delta(\varepsilon)>0$ such that
\begin{equation}
\label{eqn:Burgess} \left|\sum_{x_1\le N_1}\chi(x_1)\right|\le
N_1m^{-\delta}.
\end{equation}

Let ${\cal H}$ be the set of all elements of $\mathbb{Z}_m^*$ such
that for each $h\in {\cal H}$ the congruence
$$
h\equiv ux_1x_2x_3\pmod m, \quad u\in [1, U],\quad x_j\in [1, N_j]
$$
is not solvable. Therefore, since $(h,m)=1,$ we have
$$
\sum_{\chi}\sum_{u\le U}\sum_{x_1\le N_1}\sum_{x_2\le N_2}
\sum_{x_3\le N_3}\sum_{h\in {\cal
H}}\chi(ux_1x_2x_3)\overline{\chi}(h)=0.
$$
Separating the term corresponding to the principal character
$\chi=\chi_0$, we get that
\begin{equation}
\label{eqn:MainTerm} UN_1N_2N_3|{\cal H}|\lessapprox
\sum_{\chi\not=\chi_0}\left|\sum_{u\le
U}\chi(u)\right|\left|\sum_{x_1,
x_2,x_3}\chi(x_1x_2x_3)\right|\left|\sum_{h\in {\cal
H}}\chi(h)\right|.
\end{equation}
Here we used the fact that the intervals $[1,U]$ and $[1,N_j]$
contain accordingly $Um^{o(1)}$ and $N_j m^{o(1)}$ numbers coprime
to $m$ (consider, for example, the primes of these intervals that
are not divisors of $m$).

The set of nonprincipal characters $\chi\pmod m$ we split into two
subsets:
\begin{eqnarray*}
{\cal A}:=\{\,\chi\pmod m: \, \left|\sum_{u\le
U}\chi(u)\right|\ge Um^{-\delta/4n}\,\},\\
{\cal B}:=\{\,\chi\pmod m: \, \left|\sum_{u\le U}\chi(u)\right|<
Um^{-\delta/4n}\,\}.
\end{eqnarray*}
It follows that
\begin{equation} \label{eqn:Afirst} \frac{|{\cal
A}|U^{2n}m^{-\delta/2}}{\varphi(m)}\le
\frac{1}{\varphi(m)}\sum_{\chi}\left|\sum_{u\le
U}\chi(u)\right|^{2n}.
\end{equation}
The right hand side of this inequality is not greater than the
number of solutions of the congruence
$$
u_1u_2\cdots u_n\equiv u_{n+1}u_{n+2}\cdots u_{2n}\pmod m,\quad
u_j\in [1, U].
$$
In view of $U^n\le m,$ this congruence implies the equality
$$
u_1u_2\cdots u_n=u_{n+1}u_{n+2}\cdots u_{2n}.
$$
Since any positive integer $x$ has  $x^{o(1)}$ divisors, the number
of solutions of this equation is $U^{n+o(1)}.$ Thus,
from~\eqref{eqn:Afirst} it follows that
$$
\frac{|{\cal  A}|U^{2n}m^{-\delta/2}}{m}\lessapprox U^{n}.
$$
Since $U^n\approx m,$ we get that
$$
|{\cal A}|\lessapprox m^{\delta/2}.
$$
Therefore, applying Burgess bound to the sum over $x_1,$ we obtain
that
\begin{eqnarray*}
&& \sum_{\chi\in {\cal A}}\left|\sum_{u\le
U}\chi(u)\right|\left|\sum_{x_1\le
N_1}\chi(x_1)\right|\left|\sum_{x_2\le
N_2}\chi(x_2)\right|\left|\sum_{x_3\le
N_3}\chi(x_3)\right|\left|\sum_
{h\in {\cal H}}\chi(h)\right|\lessapprox \\
&&\qquad \qquad \lessapprox m^{\delta/2}UN_1m^{-\delta}N_2 N_3|{\cal
H}|\lessapprox m^{-\delta/2}UN_1N_2N_3|{\cal H}|.
\end{eqnarray*}
Inserting this into the inequality~\eqref{eqn:MainTerm}, we see that
the sum over $\chi\in {\cal A}$ never dominates, and we therefore
get
$$
UN_1N_2N_3|{\cal H}|\lessapprox \sum_{\chi\in {\cal
B}}\left|\sum_{u\le U}\chi(u)\right|\left|\sum_{x_1, x_2,
x_3}\chi(x_1x_2x_3)\right|\left|\sum_{h\in {\cal H}}\chi(h)\right|.
$$
The sum over $u$ we estimate in accordance with the definition of
the set ${\cal B}.$ This implies, after cancelation by $U,$
$$
N_1N_2N_3|{\cal H}|\lessapprox m^{-\delta/4n}\sum_{\chi\in {\cal
B}}\left|\sum_{x_1, x_2, x_3}\chi(x_1x_2x_3)\right|\left|\sum_{h\in
{\cal H}}\chi(h)\right|.
$$
Now extending the summation over $\chi\in {\cal B}$ to the set of
all characters $\chi\pmod m$ and then applying the Cauchy-Schwarz
inequality, we deduce
\begin{eqnarray*}
&&\sum_{\chi\in {\cal B}}\left|\sum_{x_1, x_2,
x_3}\chi(x_1x_2x_3)\right|\left|\sum_{h\in {\cal
H}}\chi(h)\right|\le\\ \quad && \le
\left(\sum_{\chi}\left|\sum_{x_1, x_2,
x_3}\chi(x_1x_2x_3)\right|^2\right)^{1/2}\left(\sum_{\chi}\left|\sum_{h\in
{\cal H}}\chi(h)\right|^2\right)^{1/2}\le \sqrt{mIm|{\cal H}|},
\end{eqnarray*}
where $I$ is the number of solutions of the congruence
\begin{equation}
\label{eqn:congrI} x_1x_2x_3\equiv y_1y_2y_3\pmod m,\quad x_j,y_j\in
[1,N_j].
\end{equation}
Thus,
\begin{equation}
\label{eqn:NNNH} N_1N_2N_3|{\cal H}|\lessapprox
m^{-\delta/4n}\sqrt{mIm|{\cal H}|}.
\end{equation}
Now we write the congruence~\eqref{eqn:congrI} as the equation
$$
x_1x_2x_3 = y_1y_2y_3+mt, \quad x_j,y_j\in [1,N_j],\quad |t|\le
N_1N_2N_3/m
$$
and observe that if we fix the quadruple $(y_1, y_2, y_3, t)$ with
$t\ge 0,$ then this equation will have $m^{o(1)}$ solutions in
variables $x_1,x_2,x_3.$ Since $N_1N_2N_3>m^{1+\varepsilon},$ there
are less than $2(N_1N_2N_3)^2m^{-1}$  collections of such
quadruples. Therefore,
$$
I\lessapprox (N_1N_2N_3)^2m^{-1}.
$$
Plugging this into~\eqref{eqn:NNNH}, we obtain
$$
|{\cal H}|\lessapprox m^{-\delta/4n}\sqrt{m|{\cal H}|}.
$$
This implies $|{\cal H}|\lessapprox m^{1-\delta/2n}$ and finishes
the proof of Theorem~\ref{thm:main1}.

\section{Proof of Theorem~\ref{thm:main11}}

The proof is the same as the one of Theorem~\ref{thm:main1}, where
Lemma~\ref{lem:main1} should be replaced with the following one:
\begin{lemma}
\label{lem:main11} Let $N_1N_2N_3N_4>m^{1+\varepsilon}.$ Then there
are only $O(m^{1-\varepsilon_1})$ elements  $\lambda\in
\mathbb{Z}_m^*$ such that
$$
\lambda\not \in\{x_1x_2x_3x_4\pmod m: \quad x_j\in [1,N_j]\}.
$$
\end{lemma}
The proof of Lemma~\ref{lem:main11} follows the same lines as the
proof of Lemma~\ref{lem:main1}. Here one uses Burgess' character sum
estimate over the interval of length $N_1>m^{(1+\varepsilon)/4}$
(such an estimate is guaranteed by the fact that $m$ is cubefree).

\section{Proof of Theorem~\ref{thm:main2}}

We assume that $\varepsilon$ is as small positive quantity as we
need below. Let $J$ be the number of solutions of the congruence
$$
p_1p_2(p_3+k)\equiv \lambda\pmod p, \qquad p_1,p_2,p_3\le N.
$$
Expressing $J$ via character sum estimates and separating the
contribution from the principal character we get, for some
$\delta'>0,$ that
$$
J=\Bigl(1+O(p^{-\delta'})\Bigr)\frac{\pi(N)^3}{p}+Error,
$$
where
$$
|Error|\ll \frac{1}{p}\sum_{\chi\not=\chi_0}\Bigl|\sum_{p_1\le
N}\chi(p_1)\Bigr|^2\Bigl|\sum_{p_3\le N}\chi(p_3+k)\Bigr|.
$$
We can split the interval of summation over $p_1$ into subintervals
of the form $(N_1, N_1'],$ where $N_1<N_1'\le 2N_1<2N.$ Then
decomposing into level sets, we get
\begin{equation}
\label{eqn:Error} |Error|\ll \frac{1}{p}RV_1^2V_2(\log q)^3,
\end{equation}
where $R$ is the number of non-principal characters $\chi$ for which
$$
V_1\le \Bigl|\sum_{p_1\sim N_1}\chi(p_1)\Bigr|\le 2V_1,\quad V_2\le
\Bigl|\sum_{p_3\le N}\chi(p_3+k)\Bigr|\le 2V_2.
$$
If  $V_1\le p^{\frac{5}{16}}N^{\frac{5}{12}+0.01\varepsilon},$ then
from~\eqref{eqn:Error} we get
$$
|Error|\lessapprox
\left(\frac{RV_1^2}{p}\frac{RV_2^2}{p}\right)^{1/2}p^{\frac{5}{16}}N^{\frac{5}{12}+0.01\varepsilon}.
$$
Since
$$
RV_1^2\le \sum_{\chi}\left|\sum_{p_1\sim N_1}\chi(p_1)\right|^2\le
pN_1,\qquad RV_2^2 \le \sum_{\chi}\left|\sum_{p_3\le
N}\chi(p_3)\right|^2\le pN,
$$
we get that
$$
|Error|\lessapprox Np^{5/16}N^{5/12+0.01\varepsilon}
$$
and thus $Error = o(\pi(N)^3/p).$

If $V_1\ge p^{\frac{5}{16}}N^{\frac{5}{12}+0.01\varepsilon},$ then
in~\eqref{eqn:Error} we apply  Vinogradov's bound~\eqref{eqn:Vin} to
get
$$
|Error|\lessapprox \frac{RV_1^2}{p}p^{1/4}N^{2/3}.
$$
Then we use the large values estimate~\eqref{eqn:Hux1} to bound
$RV_1^2:$
$$
RV_1^2\lessapprox N^2+\frac{pN^4}{V^4}\lessapprox
\frac{pN^4}{p^{5/4+0.04\varepsilon}N^{5/3}}.
$$
The result now follows.

\section{Remarks}

Theorems~\ref{thm:main1},~\ref{thm:main11} can be included into a
more general statement. For instance, let $k$ be fixed,
$N_1,N_2,\ldots, N_k$ be positive integers such that
$N_1>m^{1/3+\varepsilon}$ and $N_1N_2\cdots N_k>m^{1+\varepsilon}.$
Then the set
$$
\{\,x_1x_2\cdots x_k \pmod m: \quad x_j\in [1, N_j]\,\}
$$
contains all, but $O(m^{1-\delta})$ elements of $\mathbb{Z}_m.$ In
case of cubefree $m$ the condition $N_1>m^{1/3+\varepsilon}$ can be
replaced by $N_1>m^{1/4+\varepsilon}.$

We can state Theorem~\ref{thm:main2} in the following form. Let
$0\le \alpha< 1,\, 0\le \beta<1$ be fixed nonnegative real numbers.
Define
$$
\theta=\max\Bigl\{\frac{\alpha}{1-\beta},\,\,
\frac{5+\alpha}{7-\beta}\Bigr\}.
$$
Let $p>N>p^{\theta+\varepsilon}$ and let for any nonprincipal
character $\chi\pmod p$ we have
$$
S_N\lessapprox p^{\alpha}N^{\beta}.
$$
Then the congruence
$$
p_1p_2(p_3+k)\equiv \lambda\pmod p
$$
has $(1+o(1))\pi(N)^3/p$ solutions in primes $p_1, p_2, p_3\le N.$
The proof is the same as the proof of Theorem~\ref{thm:main2} (one
considers the cases $V_1\le
p^{\frac{1+\alpha}{4}}N^{\frac{1+\beta}{4}+0.01\varepsilon}$ and
$V_1\ge
p^{\frac{1+\alpha}{4}}N^{\frac{1+\beta}{4}+0.01\varepsilon}$). In
view of~\eqref{eqn:Vin} the pair $(\alpha,\beta)=(1/4, 2/3)$ is
acceptable, which produces $\theta=63/76.$ It would be interesting
to obtain  pairs $(\alpha,\beta)$ which would improve our exponent
$63/76.$

\end{document}